\documentclass{article}
\usepackage{amssymb}
\usepackage{amsmath}

\newtheorem{theorem}{Theorem}[section]

\newtheorem{proposition}[theorem]{Proposition}
\newtheorem{corollary}[theorem]{Corollary}

\newtheorem{preexample}{Example}[section]
\newenvironment{example}{\begin{preexample}}{\end{preexample}}
\newtheorem{preremark}{Remark}

\newenvironment{dem}
{{\bf Proof:}}
{\qquad \hspace*{\fill} $\Box$}%

\newcommand{\fg}{\mathfrak{g}}%
\newcommand{\ad}{\operatorname{ad}}%

\newcommand{\inner}{\operatorname{int}}%
\newcommand{\SC}{\mathcal{S}}%
\newcommand{\UC}{\mathcal{U}}%
\newcommand{\DC}{\mathcal{D}}%
\newcommand{\XC}{\mathcal{X}}%
\newcommand{\YC}{\mathcal{Y}}%
\newcommand{\AC}{\mathcal{A}}%
\newcommand{\R}{\mathbb{R}}%

\title{The General Solution for Affine Control Systems on Lie Groups}

\author{Jo\~ao Paulo Lima de Oliveira \and Alexandre J. Santana\\
	Departamento de Matem\'atica, Universidade Estadual de Maring\'a\\
	Maring\'a, Brazil \and Sim\~ao N. Stelmastchuk\\
	Universidade Federal do Paran\'a, Jandaia do Sul, Brazil}

\begin{document}

\maketitle

\begin{abstract}
	The purpose of this paper is to present explicitly the solution curve for affine control systems on Lie groups under the assumption that automorphisms associated to the linear vector fields commutes. If we assume that the derivations associated to linear vector fields are inner, we obtain a simpler solution and we show some results of controllability. To end, we work with conjugation by homomorphism of Lie groups between affine systems. 
\end{abstract}
	
\section{Introduction}

Let $G$ be a connected Lie group and denote by $\fg$  the Lie algebra of the right invariant vector fields of $G$. A vector field is called {\it linear} if its flow $\varphi_t$ is a one parameter group of automorphism of $G$. An affine system on $G$ is a control system of the form
\[
	\Sigma\colon\displaystyle\frac{dg}{dt}=({\cal{X}}+Y)(g)+\displaystyle\sum_{j=1}^{m} u_{j}({\cal{X}}_j+Y_j)(g)
\]
where ${\cal X}, {\cal X}_1,\ldots, {\cal X}_m$ are linear vector fields, $Y, Y_1,\ldots, Y_m$ right invariant vector fields on $G$, and $u = (u_1, \ldots, u_m)$ is an admissible control. Note that an affine system on a Lie group $G$ is an extension of one on $\R^n$ (see for instance \cite{jouan2}). Controllability and conjugation of affine systems has been study for some authors, for example, Jurjdevic and Sallet in \cite{JurdSall}, Kara and San Martin in \cite{Karas}, Rocio, Santana and Verdi in \cite{Rocio}, Jouan in \cite{jouan2}, and, more recently, Ayala, Da Silva and Ferreira in \cite{Ayala}.

Our first purpose is to study the solution of $\Sigma$. To construct the solution (see Theorem \ref{Formgersol}) we use a technique developed in \cite{Cardetti} and improved in \cite{Oliveira} (see Section 2 for more details). In particular, considering $\R^n$ as Lie group it is not difficult to see that our solution is an extension of that one presented in \cite[ch.3]{Agrachev}. However, we need the assumption: $[{\cal X}_i, {\cal X}_j]=0$, for $i,j=0,1,\ldots,m$, with $\XC_0 = \XC$. We observe that this assumption is not an obstruction. In fact, in a direct product of Lie groups we can construct an affine system that satisfies this. It is the case of compact Lie groups (see for instance \cite{Arvanitoyeorgos}). Furthermore, we can find this assumption naturally in semisimple Lie groups (see for instance Theorem 3.11 in \cite{SanMartin2}).  

Knowing that each linear vector field $\XC$ yields the derivation of Lie algebra given by $\DC(Y) = - [\XC,Y]$, $Y \in \fg$. We assume that derivations $\DC_i$ yielded by linear vector fields $\XC_i$, $i= 1, \ldots, 0$, are inner. We remember that every derivation is inner. Under this assumption, we prove in Theorem \ref{solsemsimpgroup} that the solution of $\Sigma$ is written in a simpler way. This also allows us to consider an invariant system. Through the invariant system we can establish conditions to controllability of $\Sigma$, see Theorem \ref{teo1}, Corollary \ref{cor1} and Theorem \ref{teo2}. 

Finally, we establish a conjugation between affine systems. Our idea is based in a conjugation by homomorphism of linear system presented in \cite{solutions}. We show that a necessary and sufficient condition for affine system to be conjugate by homomorphism of Lie groups is that the flows of linear vector fields commute and invariant vector fields are related (as one can see in Theorem \ref{teo3}).

The paper is organized as follows, in the second section we establish some basic facts about linear vector fields, we construct the solution of affine system and give some examples. In third section, under assumption that derivations are inner, we rewrite the solution of affine system and obtain results of controllability. Finally, in the last section we study a conjugation between affine systems by homomorphisms.

\section{Solution for Affine Control Systems on Lie Groups}

In this section, we construct a solution for an affine control system on a Lie group. We begin by introducing linear vector fields on Lie groups. 

Let $G$ be a real, connected Lie group and let us denote by $\fg$ its Lie algebra. A vector field $\XC$ on $G$, whose flow is denoted by $\varphi_t$, is called {\it linear} if it satisfies some of the equivalent sentences:
\begin{description}
	\item {(i)} for all $t \in \R$, $\varphi_t$ is an automorphism of $G$;
	\item{(ii)} for all $Y\in\fg$, $[\XC,Y]\in\fg$;
	\item{(iii)} for all $g,h \in G$, $\XC(gh) = d(R_h)_g\XC(g) + d(L_g)_h\XC(h)$.
\end{description}
Note that $\XC(e) = 0$, where $e$ is the identity of $G$ and any linear vector field $\XC$ define a derivation $\DC \colon\mathfrak{g}\to\mathfrak{g}$ by $\DC =-\ad({\XC})$. This derivation satisfies the condition $d(\varphi_t)_e=e^{t\DC}$, for all $t\in\R$.

Our next step is to define an affine control system on $G$. Let $\XC, \XC_1,\ldots, \XC_m$ be linear vector fields and $Y, Y_1,\ldots, Y_m$ right invariant vector fields on $G$. An affine control system, or shortly affine system, is a control system given by a family of differential equations
\begin{equation}\label{affinesystem}
	\Sigma\colon\displaystyle\frac{dg}{dt}=(\XC+Y)(g)+\displaystyle\sum_{j=1}^{m} u_{j}(\XC_j+Y_j)(g),
\end{equation}
where control functions $u\colon\R\to U\subset\R^m$ belong to a subset ${\cal U}\subset L_{loc}^{\infty}(\R;\R^m)$ of the space of the locally integrable functions. 
Let us denote by $\phi_t(g,u)$ the solution of affine system (\ref{affinesystem}) starting at $g$. The reachable set of affine system from a point $g \in G$  at time $t>0$ is defined by 
\[
	{\cal A}_t(g) = \{ h \in G: \phi_t(g,u) \}.
\]
Also, the reachable set of affine system (\ref{affinesystem}) from a point $g \in G$ is given by  ${\cal A}(g) = \bigcup_t {\cal A}_t(g)$. 

In this context we have the following particular control systems:
\begin{enumerate}
	\item Invariant system if $\XC = \XC_1 = \ldots = \XC_m = 0$;
	\item Bilinear system if $Y = Y_1 = \ldots = Y_m =0 $;
	\item Linear system if $Y = 0 $ and $\XC_1 = \ldots = \XC_m = 0$.
\end{enumerate}

Our next step is to introduce an approach to construct a solution to affine system. Our idea is to follow a technique due to Cardetti and Mittenhuber in \cite{Cardetti}. Their idea, to study local controllability for linear systems, is: first, to construct a semidirect product $\R \times_{\varphi_t} G$, where $\varphi_t$ is a linear flow; second, to lift the linear system to an invariant system on $\R \times_{\varphi_t} G$; third, to study the controllability of invariant system on $\R \times_{\varphi_t} G$ and to construct a way to project the controllability results on $G$. We want to follow this idea with a little change. Instead of consider left invariant systems we consider right invariant systems. This allows us to define a natural projection of $\R \times_{\varphi_t} G$ over $G$.

Denote by $\varphi_t, \varphi^1_t, \ldots, \varphi^m_t$ linear flows associated to the linear vector fields $\XC, \XC_1, \ldots, \XC_m$, respectively. Consider the application $\rho\colon\R^{m+1}\rightarrow Aut(G)$ defined as
\[
	\rho(t_0,t_1,\ldots,t_m)=\varphi_{t_0}\circ\varphi^1_{t_1}\circ\cdots\circ\varphi^m_{t_m}.
\]
Assuming that 
\begin{equation}\label{assumption1}
	[\XC, \XC_j]=0\ \ \mbox{and}\ \ [\XC_i, \XC_j]=0,\ \ \mbox{for}\ \ i,j = 1, \ldots, m,
\end{equation}
we have that $\rho$ is a representation of $\R^{m+1}$ into $G$. In fact,
\begin{eqnarray*}
	\rho(t_0+s_0,t_1+s_1,\ldots,t_m+s_m)&=&\varphi_{t_0+s_0}\circ\varphi^1_{t_{1}+s_{1}}\circ\cdots\circ\varphi^m_{t_{m}+s_{m}}\\
	&=&\varphi_{t_0}\circ\varphi_{s_0}\circ\varphi^1_{t_{1}}\circ\varphi^1_{s_{1}}\circ\cdots\circ\varphi^m_{t_{m}}\circ\varphi^m_{s_{m}}\\
	&=&\varphi_{t_0}\circ\varphi^1_{t_{1}}\circ\cdots\circ\varphi^m_{t_{m}}\circ\varphi_{s_0}\circ\varphi^1_{s_1}\circ\cdots\circ\varphi^m_{s_m}\\
	&=&\rho(t_0,t_1,\ldots,t_m)\circ\rho(s_0,s_1,\ldots,s_m).
\end{eqnarray*}

Assumption (\ref{assumption1}) is not rare. For example, if $G$ is a direct product $G_0 \times G_1 \times \ldots \times G_m$, then taking $\XC \in G_0$ and $\XC_i \in G_i$, for $i = 1, \ldots, m$ we can view that condition (\ref{assumption1}) is satisfied. An especial case of this is when $G$ is a compact Lie group because it is isomorphic to a direct product of simple, compact, connected and simply connected Lie groups. See for instance \cite{Arvanitoyeorgos} to view the list of simple, compact, connected and simply connected Lie groups.

We then define the semi-direct product $G\times_{\rho}\R^{m+1}$, that is, the cartesian product of $G$ and $\R^{m+1}$ endowed with the product $(g,t)(h,s)=(g\rho_{t}(h),t+s)$. This set is a Lie group and the correspondent Lie algebra is the semi-direct product of algebras $\mathfrak{g}\times_{\sigma}\R^{m+1}$, where $\sigma\colon\R^{m+1}\rightarrow Der(\mathfrak{g})$ is defined as
\[
	\sigma_t(Y)=\sigma(t)(Y)=ad_{\sum t_{i}{\cal X}_i}(Y),
\]
for $t=(t_0,\ldots,t_m)$ and $Y\in\fg$. Our idea is to construct an invariant control system on $G\times_{\rho}\R^{m+1}$ from affine system (\ref{affinesystem}), but before that we need to establish some results about invariant vector fields of the Lie algebra $\mathfrak{g}\times_{\sigma}\R^{m+1}$.
\begin{proposition}
	Let $(Y,t_0,\ldots,t_m)$, $(W,s_0,\ldots,s_m)$ be vector fields in $\mathfrak{g}\times_{\sigma}\R^{m+1}$. Then
	\[
		\left[\left(Y,t_0,\ldots,t_m)(W,s_0,\ldots,s_m\right)\right]=\left([Y+\displaystyle\sum_{i=0}^m t_{i}{\cal X}_i,W+\displaystyle\sum_{j=0}^m s_{j}{\cal X}_j],0,\ldots,0\right).
	\]
\end{proposition}
\begin{dem}
	We begin by computing 
	\begin{eqnarray*}
		\left[(Y,t)(W,s)\right]&=&([Y,W]+\sigma_t(W)-\sigma_s(Y),[t,s])\\
		&=&\left([Y,W]+ad_{\sum t_{i}{{\cal X}_i}}(W)-ad_{\sum s_{j}{{\cal X}_j}}(Y),0\right).
	\end{eqnarray*}
	Adding $\left[\sum t_i{\cal X}_i,\sum s_j{\cal X}_j\right]=0$ in the first coordinate we get
	\[
		\left[(Y,t)(W,s)\right] = \left([Y+\displaystyle\sum t_{i}{\cal X}_i,W+\displaystyle\sum s_j{\cal X}_j],0,\ldots,0\right).
	\]
\end{dem}

Let $(W,0,\ldots,0)$, $(0,s_0,s_1,\ldots,s_m)\in\mathfrak{g}\times_{\sigma}\R^{m+1}$. A direct calculation proves that their exponentials are $(\exp(tW),0,\ldots,0)$ and $(e,s_0t,s_1t,\ldots,s_mt)$, respectively. This fact allows us to obtain the exponential for any invariant vector field on $G\times_{\rho}\R^{m+1}$.

\begin{proposition}\label{FormExp}
	If $(W,s_0,s_1,\ldots,s_m)$ is a vector field in $\mathfrak{g}\times_{\sigma}\R^{m+1}$, then 
	\begin{eqnarray}
		\exp(t(W,s))=\left(\displaystyle\lim_{n\rightarrow\infty}\displaystyle\prod_{i=0}^{n-1}\rho(is_0t/n,\ldots,is_mt/n)(\exp(t/n\cdot W)),s_0t,\ldots,s_mt)\right).
	\end{eqnarray} 
\end{proposition}
\begin{dem}
	We first write $(W,s_0,\ldots,s_m)=(W,0,\ldots,0)+(0,s_0,\ldots,s_m)$. Now, applying the Lie product formula we obtain
	\begin{eqnarray*}
		\exp(t(W,s))&=&\displaystyle\lim_{n\rightarrow\infty}\left(\exp(t/n\cdot W,0)\cdot \exp(0,ts/n)\right)^n\\
		&=&\displaystyle\lim_{n\rightarrow\infty}\left((\exp(t/n\cdot W),0,\ldots,0)(e,s_0t/n,s_1t/n,\ldots,s_mt/n)\right)^n\\
		&=&\displaystyle\lim_{n\rightarrow\infty}(\exp(t/n\cdot W),s_0t/n,s_1t/n,\ldots,s_mt/n)^n\\
		&=&\left(\displaystyle\lim_{n\rightarrow\infty}\displaystyle\prod_{i=0}^{n-1}\rho(is_0t/n,\ldots,is_mt/n)(\exp(t/n\cdot W)),s_0t,\ldots,s_mt)\right).
	\end{eqnarray*}
\end{dem}

Denoting by $D_{{\cal X}_0},\,\ldots,\,D_{{\cal X}_m}$ the derivations of linear vector fields ${\cal X}_0,\,\ldots,{\cal X}_m$, respectively, we can rewrite the above result as 
\begin{equation}\label{altsolution}
	\exp(t(W,s))=\left(\displaystyle\lim_{n\rightarrow\infty}\displaystyle\prod_{i=0}^{n-1}\exp\left(t/n\cdot e^{D_t} W\right),s_0t,\ldots,s_mt\right)
\end{equation}
where $D_t=\dfrac{it}{n}D_{{\cal X}_0}+\dfrac{iu_1t}{n}D_{{\cal X}_1}+\cdots+\dfrac{iu_mt}{n}D_{{\cal X}_m}$.

Our next step is describe how an invariant vector field $(W,s_0,\ldots,s_m)$ acts on an arbitrary point $(g,r_0,\dots,r_m)$.
\begin{proposition}\label{Efcamp}
	If $(W,s)\in\mathfrak{g}\times_\sigma\R^{m+1}$ and $(g,r)\in G\times_\rho\R^{m+1}$, then 
	\[
		(W,s)(g,r)=\left(W(g)+\displaystyle\sum_{i=0}^{m}s_i{\cal X}_i(g),s_0,\ldots,s_m\right),
	\]
	where $s=(s_0,\ldots,s_m)$ and $r=(r_0,\ldots,r_m)$.
\end{proposition}
\begin{dem}
	Using the right invariance, we have 
	\[
		\!\!\!\!\!\!\!\! (W,s)(g,r) = d(R_{(g,r)})(W,s) = d(R_{(g,r)})(W,0)+\displaystyle\sum_{i=0}^md(R_{(g,r)})(0,\ldots,s_i,\ldots,0).
	\]
	By definition of exponential on $G\times_\rho\R^{m+1}$, 
	\begin{eqnarray*}
		(W,s)(g,r) & = & \left.\frac{d}{dt}\left((\exp(tW),0)(g,r)\right)\right|_{t=0}+\displaystyle\sum_{i=0}^m\left.\frac{d}{dt}\left(e,0,\ldots,s_it,\ldots,0)(g,r)\right)\right|_{t=0}\\
			&=&\left.\frac{d}{dt}\left(\exp(tW)g,r\right)\right|_{t=0}+\displaystyle\sum_{i=0}^m\left.\frac{d}{dt}\left(\varphi^i_{s_it}(g),0,\ldots,s_it+r_i,\ldots,0\right)\right|_{t=0}.
	\end{eqnarray*}
	Differentiating  each term of right side with respect $t$ yields
	\[
		(W,s)(g,r) = (W(g),0)+\displaystyle\sum_{i=0}^m\left(s_i{\cal X}_i(g),0,\ldots,s_i,\ldots,0\right)=\left(W(g)+\displaystyle\sum_{i=0}^ms_i{\cal X}_i(g),s\right).
	\]
\end{dem}

Consider $\bar{\XC_j}=(0,\ldots,1,\ldots,0)$ and $\bar{Y_j}=(Y_j,0,\ldots,0)\in\fg\times_{\sigma}\R^{m+1}$, for $j=0,1,\ldots,m$, where $1$ is placed at the $j$-th position. From the previous proposition we see that, in coordinates, these fields may still be expressed as $\bar{{\cal X}_j}(g,r)=({\cal X}_j(g),0,\ldots,1,\ldots,0)$ and $\bar{Y}_j(g,r)=(Y_j(g),0,\ldots,0)$, $j=0,\ldots,m$. Hence we have (\ref{affinesystem}) the following invariant control system on $G\times_{\rho}\R^{m+1}$, associated to the affine system \ref{affinesystem}:
\[
	\bar{\Sigma}\colon\displaystyle\frac{d(g,r)}{dt}=(\bar{{\cal X}}+\bar{Y})(g,r)+\displaystyle\sum_{j=1}^{m} u_{j}(\bar{{\cal X}_j}+\bar{Y_j})(g,r).
\]
In coordinates, we have
\[
	\left(\begin{array}{c}
		dg/dt\\
		dr_0/dt\\
		dr_1/dt\\
		\vdots\\
		dr_m/dt
	\end{array}\right)=
	\left(\begin{array}{c}
		({\cal{X}}+Y)(g)+\displaystyle\sum_{j=1}^{m} u_{j}({\cal X}_j+Y_j)(g)\\
		1\\
		u_1\\
		\vdots\\
		u_m
	\end{array}\right).
\]

This means that the invariant control system $\bar{\Sigma}$ was built to satisfy $\pi(\bar{\Sigma})=\Sigma$, where $\pi\colon G\times_{\rho}\R^{m+1}\to G$ is the projection on the first coordinate. If we denote $\bar{{\cal A}_t}(g,r)$ the reachable set of a point $(g,r)$ in time $t>0$ to the invariant system $\bar{\Sigma}$, then $\pi(\bar{{\cal A}_t}(g,r))={\cal A}_t(g)$. 


We now are in position to prove our main result.

\begin{theorem}\label{Formgersol}
	Consider the curve
	\begin{equation}\label{SolAffSystem}
		\phi_t(e,u)=\displaystyle\lim_{n\rightarrow\infty}\displaystyle\prod_{i=0}^{n-1}\rho(it/n,iu_1t/n,\ldots,iu_mt/n)\exp\left(\dfrac{t}{n}\displaystyle\sum_{j=1}^{m}u_jY_j\right),
	\end{equation}
	where $u=(1, u_1,\ldots,u_m) \in \R^{m+1}$. Then $\phi_t(e,u)$ is the solution of the dynamical system
	\begin{equation}\label{eq1}
		\displaystyle\frac{dg}{dt}=({\cal{X}}+Y)(g)+\displaystyle\sum_{j=1}^{m} u_{j}({\cal{X}}_j+Y_j)(g)
	\end{equation}
	with initial condition $\phi_0(e,u) = e$.
\end{theorem}
\begin{dem}
	We begin by writing $W=Y+\displaystyle\sum_{j=1}^{m}u_jY_j$. From Proposition \ref{FormExp} we see that $\exp(t(W,u))=(\phi_t(e,u),t,u_1t,\ldots,u_mt)$. Now Proposition \ref{Efcamp} leads to \\
	\\
	$ (W,u)(\phi_t(e,u),t,u_1t,\ldots,u_mt) = $
	\begin{eqnarray*}
		&=& (W(\phi_t(e,u))+{\cal X}(\phi_t(e,u))+\displaystyle\sum_{j=1}^mu_j{\cal X}_j(\phi_t(e,u)),1,u_1,\ldots,u_n)\\
		&=& \left(({\cal{X}}+Y)(\phi_t(e,u))+\displaystyle\sum_{j=1}^{m} u_{j}({\cal{X}}_j+Y_j)(\phi_t(e,u)),\ldots,u_m\right).
	\end{eqnarray*}
	On the other hand, 
	\[
		(W,1,\ldots,u_m)(\phi_t(e,u),t,u_1t,\ldots,u_mt)=(d\phi_t(e,u)/dt,1,\ldots,u_m).
	\]
	So, in coordinates, it follows
	\[
		\left(\begin{array}{c}
			d\phi_t(e,u)/dt\\
			1\\
			u_1\\
			\vdots\\
			u_m
		\end{array}\right)=
		\left(\begin{array}{c}
			({\cal{X}}+Y)(\phi_t(e,u))+\displaystyle\sum_{j=1}^{m} u_{j}({\cal X}_j+Y_j)(\phi_t(e,u))\\
			1\\
			u_1\\
			\vdots\\
			u_m
		\end{array}\right).
	\]
	Taking the projection on the first coordinate we see that the curve $\phi_t(e,u)$ satisfies the differential equation (\ref{eq1}). Since $\phi_0(e,u)=e$, we conclude that $\phi_t(e,u)$ is the solution of the system at the identity.
\end{dem}

The theorem shows the solution of affine system at identity. However, it is possible to describe the solution of affine system in an arbitrary point $g \in G$.

\begin{corollary}\label{Solarbpoint}
	If $u=(u_1,u_2,\ldots, u_m) \in \R^{m}$ a constant admissible control, then the solution of the affine system (\ref{affinesystem}) at an arbitrary point $g\in G$ is given by 
	\[
		\phi_{t}(g,u)=\phi_t(e,u)\rho(t,u_1t,\ldots,u_mt)(g).
	\]
\end{corollary}
\begin{dem}
	Consider a point $(g,r)\in G\times_{\rho}\R^{m+1}$, where $r=(r_0,\ldots,r_m)$ is arbitrary. Let us denote by $\psi_{t}(g,r,u)$ the solution of the system $\bar{\Sigma}$. Since $\bar{\Sigma}$ is an invariant system, it follows $\psi_{t}(g,r,u)=\psi_{t}(e,0,u)(g,r)$. On the other hand, we have that $\pi(\psi_{t}(g,r,u))=\phi_{t}(g,u)$. So
	\begin{eqnarray*}
		\phi_{t}(g,u)&=&\pi\left(\psi_{t}(e,0,u)(g,r)\right)\\
		&=&\pi(\left(\phi_t(e,u),t,u_1t,\ldots,u_mt\right)(g,r_0,\ldots,r_m))\\
		&=&\pi(\phi_t(e,u)\rho(t,u_1t,\ldots,u_mt)(g),t+r_0,\ldots,u_mt+r_m)\\
		&=&\phi_t(e,u)\rho(t,u_1t,\ldots,u_mt)(g).
	\end{eqnarray*}
\end{dem}

The previous theorem and its corollary allow us to describe the solution curve for affine control systems by applying the cocycle property, as soon as admissible controls are piece-wise constant. In fact, without loss of generality, consider an admissible control $u\colon[0,t+s]\to U\subset\R^m$ with $t,s\in\R$ given by 
concatenation 
\[
	u(r) = \left\{
	\begin{array}{ll}
		u_1 &, \mbox{if } \,\, r \in [0,t]\\
		u_2 &, \mbox{if } \,\, r \in [t, t+s],
	\end{array}
	\right.
\]
where $u_1$ and $u_2$ are constants. By cocycle property, $\phi_{t+s}(g,u)=\phi_s(\phi_t(g,u_1),u_2)$. Applying Corollary \ref{Solarbpoint} yields
\[
	\phi_{t+s}(g,u)=\phi_s(e,u_2)\rho(t(1,u_2))(\phi_t(g,u_1)).
\]
This process can be extended for any admissible control since they can be written as a concatenation of an arbitrary quantity of constant controls. 

To end this section, we presents some examples.

\begin{example}[{\bf Invariant Control Systems}]
	An invariant control system is given by 
	\[
		\displaystyle\frac{dg}{dt}=Y(g)+\displaystyle\sum_{j=1}^{m} u_{j}Y_j(g),
	\]
	where $Y, Y_1, \ldots, Y_m$ are right invariant vector fields on $G$ and $u =(u_1, \ldots, u_m)$ is an admissible control. It is clear that it is a particular case of an affine system. In particular, we can assume that the representation $\rho$ is the identity map. Using Theorem \ref{Formgersol} we recover the well known solution
	\[
		\phi_t(e,u) = \displaystyle\lim_{n\rightarrow\infty}\displaystyle\prod_{i=0}^{n-1}\exp\left(\dfrac{t}{n}\displaystyle\sum_{j=1}^{m}u_jY_j\right)
		 =  \exp\left(t\displaystyle\sum_{j=1}^{m}u_jY_j\right).
	\]
\end{example}

\begin{example}[{\bf Bilinear Control Systems}]\label{solbilsys}
	A bilinear control system is a control system defined by 
	\[
		\displaystyle\frac{dg}{dt}={\cal X}(g)+\displaystyle\sum_{j=1}^{m} u_{j}{\cal X}_j(g),
	\]
	where $\XC, \XC_1, \ldots, \XC_m$ are linear vector fields on $G$ and $u =(u_1, \ldots, u_m)$ is an admissible control. Since the identity is a singularity point, we describe the solution at an arbitrary point $g\in G$. From Corollary \ref{Solarbpoint} it follows immediately that
	\[
		\phi_{t}(g,u) = \rho(t,u_1t,\ldots,u_mt)(g).
	\]
	In particular, if we consider a bilinear control system on $\R^n$ given by
	\[
		\dfrac{dx}{dt}=\left(A+\displaystyle\sum_{i=1}^mu_iB_i\right)x.
	\]
	It follows that the solution at a point $x$ is written as 
	\[
		\phi_{t}(x,u) = e^{tA}e^{u_1tB_1}\cdots e^{u_mtB_m}x.
	\]
\end{example}

\begin{example}[Linear control system on $Gl(n;\R)^+$]\label{exlinsys}
	Let $Gl(n;\R)^+$ be the set of all $n\times n$ real matrices with positive determinant and $\mathfrak{gl}(n;\R)$ its Lie algebra. For $A \in \mathfrak{gl}(n;\R)$ the vector field $\XC_A(g)=Ag-gA$ is linear, and its linear flow is given by $\varphi_t(g)=e^{tA}\cdot g\cdot e^{-tA}$. Consider $B_1,\ldots,B_m\in \mathfrak{gl}(n;\R)$, then they are right invariant vector fields defined by $B_j(g)=B_jg$. Define a linear control system on $Gl(n;\R)^+$ by 
	\begin{equation}\label{linearonGL}
		\dfrac{dg}{dt}= \XC_A(g)+\displaystyle\sum_{j=1}^mu_jB_j(g).
	\end{equation}
	We want to apply Theorem \ref{Formgersol} to find the solution of the above linear control system. Note that $\rho(t,t_1,\ldots, t_m) = \varphi_t$. Thus
	\[
		\phi_t(e,u)=\displaystyle\lim_{n\rightarrow\infty}\displaystyle\prod_{i=0}^{n-1}\varphi_{it/n}\exp\left(\dfrac{t}{n}\displaystyle\sum_{j=1}^{m}u_jY_j\right).
	\]
	Then 
	\[
		\phi_t(e,u) = e^{t\left(A+\sum u_jBj\right)}e^{-tA}.
	\]
\end{example}

\section{Inner Derivation Case}

In this section we study the solution and controllability of affine systems when derivations associated to linear vector fields are inner. Let $\XC_0=\XC, \XC_1, \ldots, \XC_m$ be linear vector fields on $G$. Under our assumption, for each $i=0,\ldots,m$, there is a right invariant vector field $X_i\in\fg$ such that $\DC_i= \ad(X_i)$, where $\DC_i$ is the derivation associated to $\XC_i$, respectively. This fact implies that $\XC_i=X_i+dIX_i$, where $dIX_i$ is the left invariant vector field induced by $I\colon G\to G$, $I(g)=g^{-1}$. As a particular case, if $G$ is a semisimple Lie group, every derivation is inner (see for instance \cite{SanMartin3}).

we improve the description of the solution of an affine system, under our assumption. After, we relate it to the solution of an associated invariant system. Begin recalling that each linear flow $\varphi_t^i$ can be written as $\varphi_t^i(g) = \exp(tX_i)g\exp(-tX_i)$ for $i=1, \ldots ,m$ (see for instance \cite{Ayala} or \cite{Jouan}), and that if
\[
	[\XC_i,\XC_j ]= 0 \mbox{\,\,\,\,\,then\,\,\,\,\,} [X_i,X_j ]= 0, \ \ \mbox{for} \,\, i,j = 0,1,\ldots, m. 
\]

\begin{theorem}\label{solsemsimpgroup} 
	Under the above assumption, the solution of the affine system (\ref{affinesystem}) is written as 
	\begin{equation}\label{solsemsimp}
		\phi_t(e,u) = \exp\left(tX+tY+\displaystyle\sum_{j=1}^mu_jt(X_j+Y_j)\right)\exp\left(\displaystyle-t\sum_{i=0}^{m}u_iX_i\right),
	\end{equation}
	where $u$ is a constant admissible control.
\end{theorem}
\begin{dem}
	We first write $W=Y+\displaystyle\sum u_jY_j$. Consider $u_0 =1$. Since $\varphi_t^i(g)=\exp(tX_i)g\exp(-tX_i)$ for each $\varphi_i$, it follows that 
	\begin{eqnarray*}
		\phi_t(e,u)&=&\displaystyle\lim_{n\rightarrow\infty}\displaystyle\prod_{i=0}^{n-1}\rho(it/n,iu_1t/n,\ldots,iu_mt/n)\exp\left(\dfrac{t}{n}W\right)\\
		&=&\displaystyle\lim_{n\rightarrow\infty}\displaystyle\prod_{i=0}^{n-1}\left(\displaystyle\prod_{k=0}^m\exp\left(\dfrac{iu_kt}{n}X_k\right)\exp\left(\dfrac{t}{n}W\right)\displaystyle\prod_{k=0}^m\exp\left(-\dfrac{iu_{m-k}t}{n}X_{m-k}\right)\right)\\
		&=&\displaystyle\lim_{n\rightarrow\infty}\displaystyle\prod_{i=0}^{n-1}\left(\exp\left(\displaystyle\sum_{k=0}^m\dfrac{iu_kt}{n}X_k\right)\exp\left(\dfrac{t}{n}W\right)\exp\left(-\displaystyle\sum_{k=0}^m \dfrac{iu_{k}t}{n}X_{k}\right)\right),
	\end{eqnarray*}
	where we use the fact that $[X_i,X_j]=0$ for $i,j =0,1, \ldots,n$. Computing the product we get \\
	\\
	$\phi_t(e,u) =$
	\[
		\displaystyle\lim_{n\rightarrow\infty}\left(\exp\left(\dfrac{t}{n}W\right)\exp\left(\displaystyle\sum_{k=0}^m\dfrac{u_kt}{n}X_k\right)\right)^{n-1}\exp\left(\dfrac{t}{n}W\right)\exp\left(\displaystyle\sum_{k=0}^m \dfrac{(1-n)u_{k}t}{n}X_{k}\right).
	\]
	Inserting
	\[
		\exp\left(\dfrac{t}{n}W\right)\exp\left(\displaystyle\sum_{k=0}^m\dfrac{u_kt}{n}X_k\right)\exp\left(\displaystyle-\sum_{k=0}^m\dfrac{u_kt}{n}X_k\right)\exp\left(-\dfrac{t}{n}W\right)
	\]
	in right side of the above equality we have
	\begin{eqnarray*}
		\phi_t(e,u) & = &\displaystyle\lim_{n\rightarrow\infty}\left(\exp\left(\dfrac{t}{n}W\right)\displaystyle\prod_{k=0}^m\exp\left(\dfrac{u_kt}{n}X_k\right)\right)^n\exp\left(-t\displaystyle\sum_{i=0}^{m}u_iX_i\right).
	\end{eqnarray*}
	Finally, applying the Lie product formula we obtain
	\[
		\phi_t(e,u)=\exp\left(tX+tY+\displaystyle\sum_{j=1}^mu_jt(X_j+Y_j)\right)\exp\left(-t\displaystyle\sum_{i=0}^{m}u_iX_i\right).
	\]
	
\end{dem}

In the remainder of this section we denote an affine system (\ref{affinesystem}) as $\Sigma_A$. From $\Sigma_A$ it is possible to yield the following right invariant control system
\[
	\Sigma_I\colon\dfrac{dg}{dt}=(X+Y)(g)+\displaystyle\sum_{j=1}^mu_j(X_j+Y_j)(g),
\]
where $X_j$ satisfies $\XC_j=X_j+dIX_j$ for $j=0,\ldots,m$. It suggests that there is a relation between affine system $\Sigma_A$ and invariant system $\Sigma_I$.

\begin{proposition}
	If $u=(u_1,\ldots, u_m)$ is a constant admissible control in an interval $[0,T]$ and if $\phi^A_t(g,u)$ and $\phi^I_t(g,u)$ are solutions of affine and invariant systems, respectively, then
	\[
		\phi^A_t(g,u)=\phi^I_t(g,u)\displaystyle\exp\left(-t\Sigma_{i=0}^mu_iX_i\right).
	\]
\end{proposition}
\begin{dem}
	We first simplify notations writing $\phi^A_t$ and $\phi^I_t$ instead of $\phi^A_t(g,u)$ and and $\phi^I_t(g,u)$, respectively. Furthermore, we write $\alpha_t = \exp\left(t\Sigma_{i=0}^m u_iX_i\right)$. Differentiating $\phi^A_t\alpha_t$ yields
	\begin{eqnarray*}
		\dfrac{d}{dt}(\phi^A_t \alpha_t) 
		& = & dR_{\alpha_t}\dfrac{d}{dt}\phi^A_t + dL_{\phi^A_t} \dfrac{d}{dt}\alpha_t\\
	  & = & dR_{\alpha_t}\left(({\cal X}+Y)\phi^A_t+\sum u_i({\cal X}_i+Y)\phi^A_t\right) + dL_{\phi^A_t}\dfrac{d}{dt}\alpha_t.
	\end{eqnarray*}
	Now, writing ${\cal X}_i=X_i+dIX_i$ and using right invariance we obtain
	\[
		\dfrac{d}{dt}(\phi^A_t \alpha_t) = \Sigma_I(\phi^I_t)+dR_{\alpha_t}\left(dIX+\sum u_idIX_i\right)\phi^A_t+ dL_{\phi^A_t}\dfrac{d}{dt}\alpha_t.
	\]
	The result follows since 
	\[
		dR_{\alpha_t}\left(dIX+\sum u_idIX_i\right)\phi^A_t + dL_{\phi^A_t}\dfrac{d}{dt}\alpha_t=0.
	\]
	The converse is proved similarly.
\end{dem}

In the following we write $\SC$ instead of $\AC$ to denote the reachable sets of the invariant system. In particular, for any $t>0$, $\SC$ is the reachable set at time $t$. The next results relates the controllability of the affine systems and the associated right-invariant ones.

\begin{theorem}\label{teo1}
	Assume that controls are piecewise constants and suppose that the right invariant system $\Sigma_I$ is controllable. The following assertions are equivalent:
	\begin{itemize}
		\item [(i)] For all control $u\in{\cal U}$ and all $t\in\R$, $\exp\left(t\left(X_0+\displaystyle\sum u_iX_i\right)\right)\in \AC$;
		\item [(ii)] $\Sigma_A$ is controllable.
	\end{itemize}
\end{theorem}
\begin{dem}
	The assertion \textit{(ii)}$\Rightarrow$\textit{(i)} is immediately. So we prove \textit{(i)}$\Rightarrow$\textit{(ii)}.\newline
	For simplicity of notation, write $X_u = X_0+\sum u_iX_i$ for an adequate control $u$. We first prove that $\Sigma_A$ is controllable from identity $e$. Given $g\in G$, there are a piecewise constant control $u$ and a time $t>0$ such that $g = \varphi^I_t(e,u) \in{\cal S}_t$. This is equivalent to $g\exp(-tX_u)\in \AC_t$. By hypothesis, $\exp(tX_u) \in \AC_s$ for some $s>0$. Thus, there exists a piecewise constant control $u'$ such that
	\begin{eqnarray*}
		g & = & g\exp(-tX_u)\rho(t,u)(\exp(tX_u)) = \phi^A_t(e,u) \rho(t,u)(\phi^A_s(e,u'))\\ 
		& = & \phi^A_t(e,u'') \rho(t,u'')(\phi^A_s(e,u'')),
	\end{eqnarray*}
	where $u''$ is the concatenation of $u$ and $u'$. From Proposition \ref{Solarbpoint} it follows that
	\[
		g = \phi^A_t(\phi^A_s(e,u''),u'') = \phi^A_{t+s}(e ,u'')\in \AC_{t+s} \subset \AC
	\]
	since $u''$ is piecewise constants. It entails that $\AC = G$, and, in consequence, $\Sigma_A$ is controllable from identity $e$. Now, we prove that $\Sigma_A$ is controllable to $e$. Set $g\in G$. By assumption, there are $t>0$ and a control $u$  such that $g^{-1} = \varphi^I_{-t}(e,u)\in{\cal S}_t$. This is equivalent to $g^{-1}\exp(-tX_u)\in{\cal A}_t$. On one hand, we have that
	\begin{eqnarray*}
		\exp(-tX_u) 
		& = & g^{-1}\exp(-tX_u)\exp(tX_u)g\exp(-tX_u)=g^{-1}\exp(-tX_u)\rho(t,u)(g)\\
		& = & \phi^A_t(e,u)\rho(t,u)(g) = \phi^A_t(g,u) \in \AC_t(g),
	\end{eqnarray*}
	where we use Proposition \ref{Solarbpoint} at last equality. On the other hand, we have $\exp(tX_u) \in \AC_s$ for some $s>0$. It means that there exists a piecewise constant control $u'$ such that $\exp(tX_{u´}) = \phi^A_t(e,u')$. Set $X_{u'}Y = X_0+\sum u'_iX_i$. Since $[X_u,X_{u'}] =0$, it follows
	\begin{eqnarray*}
		e & = & \exp(tX_u)\exp(sX_{u'})\exp(-tX_u)(\exp(-sX_{u'}))\\
		e & = & \exp(tX_u)\rho(s,u')(\exp(-tX_u))\\
		& = & \phi^A_s(e,u')\rho(s,u')(\phi^A_t(g,u)) = \phi^A_s(e,u'')\rho(s,u'')(\phi^A_t(g,u'')),
	\end{eqnarray*}
	where $u''$ is the concatenation of $u$ and $u'$.From Proposition \ref{Solarbpoint} it follows that
	\[
		e = \phi^A_t(\phi^A_s(g,u''),u'') = \phi^A_{t+s}(g ,u'')\in \AC_{t+s}(g) \subset \AC(g).
	\]
We thus conclude that $e\in{\cal A}(g)$, and the proof is complete.
\end{dem}

\begin{corollary}\label{cor1}
	Under the hypothesis of previous Theorem, if $e\in\inner(\AC_t)$ for some $t>0$, then $\Sigma_A$ is controllable.
\end{corollary}
\begin{dem}
	We write $X_u = X_0+\sum u_iX_i$ for an adequate control $u$. Let $R_u=\{t\in\R\colon\exp(tX_u)\in \AC\}$. Analysis similar to that in the proof of previous theorem shows that $R_u$ is a semigroup. By hypothesis, $\AC_t$ is a neighborhood of $e$ for some $t>0$. It implies that $0 \in S$ because $\exp(0X) = e \in \AC$. For each admissible control $u$ the curve $\exp(tX_u)$ is continuous. Then $\exp(sX_u) \in \AC_t$ for $s \in (a,b)$, where $(a,b)$ is an interval such that $0 \in (a,b)$. In particular, $(a,b) \subset A_u$. Being $S_u$ semigroup, it follows that $ S_u = \R$. It means that for all control $u\in {\cal U}$ and all $t\in\R$, $\exp\left(t\left(X_0+\displaystyle\sum u_iX_i\right)\right)\in \AC$. According to above theorem, $\Sigma_A$ is controllable.
\end{dem}

The following result generalizes, for affine control systems, Theorem 2 in \cite{Jouan}. First, we need to recall that a semigroup $S\subset G$ is said to be left reversible (resp. right reversible) if $SS^{-1}=G$ (resp. $S^{-1}S=G$). It is known that if $G$ is semi-simple with finite center, the  unique subsemigroup of $G$ with nonempty interior which is left or right reversible is $G$ itself (see for instance \cite{SanMartin1}).
\begin{theorem}\label{teo2}
	Let $G$ be a semi-simple Lie group with finite center. Suppose that $\Sigma_I$ satisfies the rank condition and $\Sigma_A$ is controllable. The following assertions are equivalent:
	\begin{itemize}
		\item [(i)] For all $u\in \UC$ and all $t\in\R$, $\exp\left(t\left(X_0+\displaystyle\sum u_iX_i\right)\right)\in \SC$;
		\item [(ii)] $\Sigma_I$ is controllable.
	\end{itemize}
\end{theorem}
\begin{dem}
	It is easy to see that that \textit{(ii)} implies \textit{(i)}. Let us prove the converse. We begin by recalling that the reachable set $\SC$ of $\Sigma_I$ is a semigroup. Now, the rank condition assures that interior of $\SC$ is non-empty interior. It is sufficient to prove that ${\cal S}$ is left reversible.
	Fix $g\in G$. There are $t>0$ and $u\in {\cal U}$ such that $g=\phi^I_t(e,u)\exp(-tX_u)$. By assumption, $\exp(-tX_u) \in \SC^{-1}$. Then $g\in {\cal S}{\cal S}^{-1}$. As $g\in G$ was chosen arbitrarily, we conclude $G\subset{\cal S}{\cal S}^{-1}$, and result follows.
\end{dem}

\section{Conjugation of affine system}

In \cite{solutions} is presented a conjugation of linear system by homomorphism. Since affine systems are a natural extension of linear system, in this section, we extent those results to affine system. 

Let $G$ and $H$ be connected Lie groups. Consider the following affine systems

\begin{equation}\label{s1}
	\displaystyle\frac{dg}{dt}=({\cal{X}}+Z)(g)+\displaystyle\sum_{j=1}^{m}u_{j}({\cal X}_j+Z_j)(g)
\end{equation}

\begin{equation}\label{s2}
	\displaystyle\frac{dh}{dt}=({\cal Y}+W)(h)+\displaystyle\sum_{j=1}^{m} u_{j}({\cal Y}_j+W_j)(h)
\end{equation}
on $G$ and $H$, respectively. Affine systems \eqref{s1}, \eqref{s2} are called conjugate if there exist a homomorphism of Lie groups $F\colon G\to H$ such that $F(\phi_t(g,u))=\theta_t(h(g),u)$, where $\phi_t(g,u)$, $\theta_t(g,u)$ are the solutions of the systems \eqref{s1}  and \eqref{s2}, respectively. 

In the following, let us denote by $\varphi^i_t$, $\psi^i_t$ the flows and by $D_{{\XC}_i}$, $D_{{\YC}_i}$ the derivations associated to the linear vector fields ${\cal X}_i$ and ${\cal Y}_i$, respectively. We give equivalent conditions for two affine systems to be conjugate. Initially, we need to extend a result of \cite{solutions}.

\begin{proposition}
	Under the above assumptions, if $F\colon G\to H$ is a homomorphism of Lie groups, then the following conditions are equivalents:
	\begin{enumerate}
		\item $F\circ\varphi^i_t=\psi^i_t\circ F$;
		\item $dF_{\varphi^i_t(g)} \XC_i(g)= \YC_i(F(g))$ for all $g \in G$;
		\item $dF_e(e^{tD_{{\XC}_i}}Z)=e^{tD_{\YC_i}}dF_eZ$.
	\end{enumerate}	
\end{proposition}
\begin{dem}
	To deduce (2) from (1), differentiate the formula in (1) with respect to $t$ to obtain $dF_{\varphi_t^i}\dfrac{d\varphi^i_t}{dt}(g)=\dfrac{d}{dt}(\psi_t^i\circ F(g))$. We thus get 
	\[
		dF_{\varphi^i_t(g)}{\cal X}_i(g)={\cal Y}_i(F(g)).
	\]
	Conversely, to deduce (1) from (2), observe that it is a direct consequence of uniqueness of solution of differential equation.
	
	Suppose now that (1) is true. Then $dF_e \circ (d\varphi_t^i)_e =d(\psi_t)_e \circ dF_e$. Since $(d\varphi_t^i)_e = e^{tD_{{\XC}_i}}$ and $(d\psi_t^i)_e = e^{tD_{{\YC}_i}}$, it follows for all $Z \in \fg$ that
	\[
		dF_e(e^{t D_{{\XC}_i}}Z)=dF_e\circ (d\varphi_t^i)_e(Z)= (d \psi_t^i)_e\circ (dF_e)(Z)=e^{t \DC_{\YC_i} } dF_e(Z).
	\]
	On the converse, to deduce (1) to (3), observe that it is a direct consequence of $d(h\circ\varphi_t^i)_e=d(\psi_t\circ h)_e$ because $G$ is connected.
\end{dem}

\begin{theorem}\label{teo3}
	Under the above assumptions, if $F\colon G\to H$ is a homomorphism of Lie groups, then following conditions are equivalents:
	\begin{enumerate}
		\item $F\left(\phi_t(g,u)\right)=\theta_t(F(g),u)$ for all $g\in G$.
		\item $F\circ\varphi^i_t=\psi^i_t\circ h$ and $dF_e Z_j(e)=W_j(e)$, for all $i,j=0,1,\ldots,m$.
	\end{enumerate}
\end{theorem}
\begin{dem}
	We first suppose that $F\left(\phi_t(g,u)\right)=\theta_t(F(g),u)$ for all $g\in G$. In particular, $F\left(\phi_t(e,u)\right)=\theta_t(e,u)$. For abbreviation, we write $\phi_t$ and $\theta_t$ instead of $\phi_t(e,u)$ and $\theta_t(e,u)$, respectively. Differentiating (1) with respect to $t$ yields
	\[
		dF_{\phi_t}\left((\XC+Z)\phi_t+\displaystyle\sum_{j=1}^{m}u_{j}(\XC_j+Z_j)\phi_t\right)=(\YC+W)\theta_t+\displaystyle\sum_{j=1}^{m} u_{j}(\YC_j+W_j)\theta_t.
	\]
	Taking $t=0$ it follows that
	\[
		dF_e\left(Z(e)+\displaystyle\sum_{j=1}^{m}u_{j}(0)Z_j(e)\right)=W(e)+\displaystyle\sum_{j=1}^{m} u_{j}(0)W_j(e).
	\]
	The above equality holds for all control $u(t)=(u_1(t),\ldots,u_m(t))$. If $u\equiv0$, then $dF_e Z(e)=W(e)$. If $u \equiv (0,\ldots,1,\ldots,0)$, then $dF_e Z_j(e)=W_j(e)$, for $j=1,\ldots,m$. Then 
	\[
		dF_{g}\left(\XC(g)+\displaystyle\sum_{j=1}^{m}u_{j}\XC_j(g)\right)=\YC(F(g))+\displaystyle\sum_{j=1}^{m} u_{j}\YC_j(F(g)).
	\]
	In the same manner we can see that $dF_{\phi_t}\XC_i(g)= \YC_i(F(g))$. From above proposition it follows that $F\circ\varphi^i_t=\psi^i_t\circ F$, for $i=0,\ldots,m$.
	
	Conversely, denote by $\rho\colon\R^{m+1}\to Aut(G)$ and $\varrho\colon\R^{m+1}\to Aut(H)$ representations associated to affine systems \eqref{s1}, \eqref{s2}, respectively. Assuming that condition (2) is true we compute
	\begin{eqnarray*}
		F(\phi_t(e,u))
		&=&F\left(\displaystyle\lim_{n\rightarrow\infty}\displaystyle\prod_{i=0}^{n-1}\rho(it/n,iu_1t/n,\ldots,iu_mt/n)\exp\left(\dfrac{t}{n}\displaystyle\sum_{j=0}^{m}u_jZ_j\right)\right)\\
		&=&\displaystyle\lim_{n\rightarrow\infty}\displaystyle\prod_{i=0}^{n-1}\varrho(it/n,iu_1t/n,\ldots,iu_mt/n)\exp\left(\dfrac{t}{n}\displaystyle\sum_{j=0}^{m}u_jdF_eZ_j\right),
	\end{eqnarray*}
	Hence $F(\phi_t(e,u))=\theta_t(e,u)$. From this last equality and Corollary \ref{Solarbpoint}  we see that
	\begin{eqnarray*}
		F(\phi_t(g,u))
		&=& F(\phi_t(e,u)\rho(t,u_1t,\ldots,u_mt)(g))\\
		&=& F(\phi_t(e,u)) F(\rho(t,u_1t,\ldots,u_mt)(g))\\
		&=& \theta_t(e,u)\varrho(t,u_1t,\ldots,u_mt)(F(g))=\theta_t(F(g),u).
	\end{eqnarray*}
\end{dem}

Next we characterize conjugation of affine system with conjugation of derivations. 

\begin{corollary}
	Under the above assumptions, if $F\colon G\to H$ is a homomorphism of Lie groups, then a necessary and sufficient condition for affine systems \eqref{s1} and \eqref{s2} to be conjugate is $dF_{e}\left(e^{D_t}Z_j\right)=e^{D'_t}W_j$, for all $j=1,\ldots,m$, $t\in\R$, where $D_t=\displaystyle\sum\dfrac{iu_kt}{n}D_{\XC_k}$ and $D'_t=\displaystyle\sum\dfrac{iu_kt}{n}D_{\YC_k}$.
\end{corollary}
\begin{dem}
	It is sufficient to show the necessary condition. Write the solution of affine system \eqref{s1}  as 
	\[
		\phi_t(e,u) = \displaystyle\lim_{n\rightarrow\infty}\displaystyle\prod_{i=0}^{n-1}\exp\left(t/n\cdot\displaystyle\sum u_je^{D_t}Z_j\right),
	\]
	by Formula \ref{altsolution}. Computing $F(\phi_t(e,u)$ yields
	\[
		F(\theta_t(e,u))=\displaystyle\lim_{n\rightarrow\infty}\displaystyle\prod_{i=0}^{n-1}\exp\left(t/n\cdot\displaystyle\sum u_jdF_{e}\left(e^{D_t}Z_j\right)\right),
	\]
	and the result follows.
\end{dem}

\begin{example}
	Consider the homomorphism $\det\colon Gl(n;\R)^+\to\R$ and the linear system \eqref{linearonGL} defined in Example \ref{exlinsys}. We construct a control system on $\R$ conjugated to it. We need a linear vector field $\YC$ and invariant vector fields $b_1,\ldots,b_m$ on $\R$ satisfying conditions
	\begin{enumerate}
		\item $\det\left(e^{tA}\cdot g\cdot e^{-tA}\right)=\psi_t(\det(g))$, for all $g\in Gl(n;\R)^+$, where $\psi_t$ is the flow of $\YC$;
		\item $d(\det)_I(B_j)=tr(B_j)=b_j$.
	\end{enumerate}
	Condition (2) gives the invariant vector fields $b_j$, $j =1, \ldots, m$. Also, condition (1) implies that $\psi_t(\det(g))=\det(g)$ for all $g$. This clearly forces $\YC = 0$. We thus conclude that the linear system \eqref{linearonGL} on $Gl(n;\R)^+$ is conjugate to the following invariant system on $\R$
	\[
		\dfrac{dx}{dt}=\displaystyle\sum_{j=1}^m tr(B_j).
	\]
	It means that linear vector field does not have importance is this conjugation. 
\end{example}

\renewcommand{\refname}{Bibliography}

\end{document}